\documentclass[11pt,a4paper]{amsart}
\usepackage{mymacros_english_paper}
\usepackage{mymacros_common}

\begin{document}
\title[A calculation of the perfectoidization of semiperfectoid rings]{A calculation of the perfectoidization of semiperfectoid rings}

\author[R. Ishizuka]{Ryo Ishizuka}
\address{Department of Mathematics, Tokyo Institute of Technology, 2-12-1 Ookayama, Meguro, Tokyo 152-8551}
\email{ishizuka.r.ac@m.titech.ac.jp}




\begin{abstract}
We show that perfectoidization can be (almost) calculated by using \(p\)-root closure in certain cases, including the semiperfectoid case.
To do this, we focus on the universality of perfectoidization and uniform completion, as well as the \(p\)-root closed property of integral perfectoid rings.
Through this calculation, we establish a connection between a classical closure operation ``\(p\)-root closure'' used by Roberts in mixed characteristic commutative algebra and a more recent concept of ``perfectoidization'' introduced by Bhatt and Scholze in their theory of prismatic cohomology.
\end{abstract}

\maketitle 

\setcounter{tocdepth}{1}
\tableofcontents
\section{Introduction}
    Let \(p\) be a prime number.
    The method of \emph{perfectoidization}, introduced by Bhatt and Scholze in \cite{bhatt2022Prismsa}, is an application of the theory of prismatic cohomology to commutative algebra.
    This yields a universal integral perfectoid ring over a ring, such as a \emph{semiperfectoid ring}, which is a derived $p$-complete ring that can be written as a quotient of an integral perfectoid ring.
    This method can be seen as a generalization of the perfect closure of positive characteristic rings. See \Cref{PerfectoidRings} for an explanation of the terminology used in perfectoid theory.

    Perfectoidization has various applications to commutative algebra, such as a new theory of almost mathematics (\cite[\S 10.1]{bhatt2022Prismsa}), a general version of the almost purity theorem (\cite[Theorem 10.9]{bhatt2022Prismsa}), much simpler proof of some previously known theorems in commutative algebra (\cite[Remark 10.13]{bhatt2022Prismsa} and \cite[Appendix A]{ma2022Analogue}), and a mixed characteristic analog of Hilbert-Kunz multiplicity and $F$-signature (\cite{cai2022Perfectoid}).

    In contrast to these applications, they are not yet widely used in commutative algebra.
    One problem is that many abstract theories, including homotopy theory, have been used, and therefore perfectoidization has a mysterious ring structure.
    To the best of the author's knowledge, perfectoidization has only been explicitly calculated in \cite[\S 2.3.1]{reinecke2020Moduli} and in the proof of \cite[Theorem 4.4]{dine2022Topologicala}.

    \subsection{\texorpdfstring{$p$}{p}-root closure}

    In this paper, we give an explicit description of the perfectoidization of semiperfectoid rings by using \(p\)-root closure.
    Before explaining our first main theorem, we recall the notion of \(p\)-root closure.

    \begin{definition}[\cite{roberts2008Root}] \label{pRootClosure}
        Let \(R\) be a \(p\)-torsion-free ring.
        We say that \(R\) is \emph{\(p\)-root closed in \(R[1/p]\)} if \(x \in R[1/p]\) satisfies \(x^{p^n} \in R\) for some \(n \geq 1\), then \(x \in R\) holds.

        The \emph{\(p\)-root closure \(C(R)\) of \(R\) in \(R[1/p]\)} is the minimal \(p\)-root closed subring of \(R[1/p]\) containing \(R\).
        Focusing on the case of ``\emph{in \(R[1/p]\)}'', Roberts provided an explicit discription of the \(p\)-root closure \(C(R)\) as follows:
        \begin{equation*}
            C(R) = \set{x \in R[1/p]}{\exists n \geq 1, x^{p^n} \in R}.
        \end{equation*}
        The term ``\emph{in \(R[1/p]\)}'' is omitted in this paper because only this case is considered.
    \end{definition}

    Here is a brief mention of the history of (\(p\)-)root closure.
    The notion of \(p\)-root closure is a special case of \emph{total \(n\)-root closure} introduced in commutative algebra by Anderson, Dobbs and Roitman \cite{anderson1990Root}.
    Previously, \emph{\(n\)-root closedness} was used, for example, by Angerm\"uller \cite{angermuller1983Root}, Anderson \cite{anderson1982Root}, Watkins \cite{watkins1982Root} and Brewer, Costa and McCrimmon \cite{brewer1979Seminormality}.
    Furthermore, its origin can be traced back to Sheldon's definition of \emph{root closedness} in \cite{sheldon1971How}.

    In the context of commutative algebra in mixed characteristic, Roberts provided the above explicit description of \(p\)-root closure and applied this even before perfectoid rings appeared.
    Most recently, (total) \(p\)-root closure has been renamed \emph{\(p\)-integral closure} by \v{C}esnavi\v{c}ius and Scholze in \cite{cesnavicius2023Purity} and is found to be more closely related to the perfectoid theory.


    \subsection{Main Theorem}

    Under this notation, our main theorems can be stated in the following forms.
    For simplicity, we only state the theorems in the case of \(p\)-torsion-free rings.
    Readers interested in the general cases may refer to the referenced statements in each theorem, in conjunction with \Cref{pTorsionCase}.

    \begin{theorem}[(\Cref{Maintheorem2}; \(p\)-torsion-free case)] \label{Maintheorem2Intro}
        Let \(R\) be a \(p\)-torsion-free ring which satisfies the following conditions:
        \begin{enumerate}
            \item The \(p\)-adic completion \(\widehat{R}\) of \(R\) has a map from some integral perfectoid ring.
            \item The perfectoidization \((\widehat{R})_{\perfd}\) of \(\widehat{R}\) is an (honest) integral perfectoid ring.\footnote{As explained in the rest of \Cref{PerfectoidRings}, the perfectoidization of an algebra over some integral perfectoid ring is only a complex. So this assumption states that \((\widehat{R})_{\perfd}\) gives an honest ring (see \Cref{PerfectoidizationUniversal}).}
            \item The \(p\)-adic completion \(\widehat{C(R)}\) of the \(p\)-root closure \(C(R)\) is an integral perfectoid ring.
        \end{enumerate}
        Then we have an isomorphism
        \begin{equation*}
            (\widehat{R})_{\perfd} \cong \widehat{C(R)}.
        \end{equation*}
    \end{theorem}

    This result clarifies the ring structure of perfectoidization by using \(p\)-root closure, which is a quite explicit closure operation.
    The left-hand side is constructed abstractly by homotopy theory, but it can be described as a \(p\)-root closure followed by a \(p\)-adic completion, which are only ring-theoretic operations.

    \subsection{Applications}

    Let \(R\) be a \(p\)-torsion-free ring such that its \(p\)-adic completion \(\widehat{R}\) is a semiperfectoid ring, that is, \(\widehat{R}\) is a quotient of some integral perfectoid ring.
    In applications of \Cref{Maintheorem2Intro}, it is crucial that \(R\) satisfies the assumptions of the theorem:
    \begin{enumerate}
        \item The semiperfectoid ring \(\widehat{R}\) has a surjective map from some integral perfectoid ring by the definition of semiperfectoid rings.
        \item The perfectoidization \((\widehat{R})_{\perfd}\) of \(\widehat{R}\) is an (honest) integral perfectoid ring by \cite[Corollary 7.3 and Proposition 8.5]{bhatt2022Prismsa} or the rest of \Cref{PerfectoidRings}.
        \item The \(p\)-adic completion \(\widehat{C(R)}\) of the \(p\)-root closure \(C(R)\) is an integral perfectoid ring by virtue of \cite[Proposition 2.1.8]{cesnavicius2023Purity} (see \Cref{pRootClosurePerfectoid}).
    \end{enumerate}
    So we can provide an explicit description of the perfectoidization of such \(R\) as follows.

    \begin{theorem}[(\Cref{MainTheorem3}; \(p\)-torsion-free case)] \label{MainTheorem3Intro}
        Let \(R\) be a \(p\)-torsion-free ring such that its \(p\)-adic completion \(\widehat{R}\) becomes a semiperfectoid ring.
        Then we have an isomorphism \((\widehat{R})_{\perfd} \cong \widehat{C(R)}\).
        In particular, any \(p\)-torsion-free and (classically) \(p\)-adically complete semiperfectoid ring \(S\) has an isomorphism \(S_{\perfd} \cong \widehat{C(S)}\).
    \end{theorem}

    In the study of commutative rings in mixed characteristic, semiperfectoid rings often appear as in \cite{ishizuka2023Mixed}.
    An application of \Cref{MainTheorem3Intro} is as follows.

    \begin{theorem}[(\Cref{ExplicitConst} and \Cref{ExplicitPerfectoidization})]
        Let \((R_0,\mfrakm,k)\) be a complete Noetherian local domain of mixed characteristic \((0, p)\) with perfect residue field \(k\) and let \(p, x_2, \dots, x_n\) be any system of generators of the maximal ideal \(\mfrakm\) such that \(p, x_2, \dots , x_d\) forms a system of parameters of \(R_0\).
        Choose compatible sequences of \(p\)-power roots 
        \begin{equation*}
            \{p^{1/p^j}\}_{j \ge 0},\{x_2^{1/p^j}\}_{j \ge 0},\ldots,\{x_n^{1/p^j}\}_{j \ge 0}
        \end{equation*}
        inside the absolute integral closure \(R_0^+\), the integral closure of \(R_0\) in the algebraic closure of the fraction field of \(R_0\).
        Set a subring \(R_{\infty}\) of \(R_0^+\) as
        \begin{equation*}
            R_{\infty} := \bigcup_{j \ge 0} R_0[p^{1/p^j},x_2^{1/p^j},\ldots,x_n^{1/p^j}] \subseteq R_0^+.
        \end{equation*}
        Then the \(p\)-adic completion \(\widehat{R}_{\infty}\) of \(R_{\infty}\) becomes a \(p\)-torsion-free semiperfectoid ring by Cohen's structure theorem (see \Cref{ExplicitConst} for details).
        Then \((\widehat{R}_\infty)_{\perfd}\) is isomorphic to the \(p\)-adic completion \(\widehat{C(R_\infty)}\) of the \(p\)-root closure \(C(R_\infty)\).
    \end{theorem}

    \subsection{Strategy of Proof}

    Let us comment on the strategy of the proof of \Cref{Maintheorem2Intro}.
    Our proof is attributed to some universalities and the following principle of ``rigidity lemma'' that makes sense in certain situations (for example \Cref{IdentityMap}).

    \begin{Lemma} \label{RigidityLemma}
        Let \(\map{f}{R}{R}\) be an endomorphism of a ring \(R\).
        Assume that \(R\) has a ``\emph{good}'' map \(S \to R\) from some ring \(S\).
        Then if \(f\) is a map of \(S\)-algebras, \(f\) is exactly the identity map. 
    \end{Lemma}
    In the proof of \Cref{Maintheorem2Intro}, two maps of rings \(\map{\varphi}{\widehat{C(R)}}{(\widehat{R})_{\perfd}}\) and \(\map{\psi}{(\widehat{R})_{\perfd}}{\widehat{C(R)}}\) are obtained through the universality of uniform completions and perfectoidizations, respectively.
    By using the lemma mentioned above, we can show that \(\psi \circ \varphi\) and \(\varphi \circ \psi\) are identity maps.
    The former is a consequence of the universality of uniform completions (\Cref{UniversalityUniformCompletion}), while the latter is a consequence of the use of ``almost elements'', which is a concept from almost mathematics (\Cref{AlmostElement} and \Cref{IdentityMap}).

    \subsection{\texorpdfstring{$p$}{p}-torsion case} \label{pTorsionCase}

    While the aforementioned theorems only deal with \(p\)-torsion-free rings, we can show that similar statements hold in general.

    For this purpose, we introduce the following symbol.
    Let \(R\) be a ring.
    The \emph{\(p\)-torsion-free quotient \(R^{\ptf}\)} is defined as the quotient ring
    \begin{equation*}
        R^{\ptf} \defeq R/R[p^\infty] \cong \Image{(R \to R[1/p])},
    \end{equation*}
    where \(R[p^\infty]\) is the ideal of all \(p^\infty\)-torsion elements of \(R\), which is the kernel of the canonical map \(R \to R[1/p]\).
    Then \(R^{\ptf}\) is a \(p\)-torsion-free ring and we have a canonical map \(R \twoheadrightarrow R^{\ptf} \hookrightarrow R[1/p]\).

    The general case of main theorems are obtained by substituting \(C(R)\) and \((\widehat{R})_{\perfd}\) for \(C(R^{\ptf})\) and \(((\widehat{R})_{\perfd})^{\ptf}\), respectively.
    Furthermore, the symbols \((-)_{\perfd}\) and \((-)^{\ptf}\) are often interchangeably because of \Cref{Interchangeable}.

    To show the general version of \Cref{Maintheorem2Intro} (i.e., \Cref{Maintheorem2}), we need to pass from the possibly \(p\)-torsion case to the \(p\)-torsion-free case as in \cite{andre2020Weak}.
    With this in mind, we show that any integral perfectoid ring can be canonically modified into a \(p\)-torsion-free integral perfectoid ring (see \Cref{UniformPerfectoid}) by using \emph{pre-perfectoid pairs} as defined in \Cref{PrePerfectoidPair}.

\subsection{Notation} \label{Notation}

    A \emph{Tate ring} is a topological ring \(A\) which has an open subring \(A_0 \subseteq A\) and an element \(t \in A_0\) such that the relative topology on \(A_0\) coincides with the \(t\)-adic topology and \(A = A_0[1/t]\) as abstract rings.
    Such an open subring \(A_0\) is called a \emph{ring of definition of \(A\)} and such an element \(t\) is called a \emph{pseudo-uniformizer of \(A\)}.
    This pair \((A_0, (t))\) of a ring and its ideal is called a \emph{pair of definition of \(A\)}.
    Note that a ring (resp., pair) of definition of \(A\) is not necessarily unique.

    Conversely, for a ring \(A_0\) and an element \(t \in A_0\), the ring \(A_0[1/t]\) becomes a Tate ring by taking \(\{t^n (A_0/A_0[t^\infty])\}_{n \geq 0}\) as a fundamental system of open neighborhoods of \(0\).
    The Tate ring \(A_0[1/t]\) has a pair of definition \((A_0/A_0[t^\infty], (t))\).
    When referring to a ring \(A_0[1/t]\) as a \emph{Tate ring}, we refer to the Tate ring that arises from the pair \((A_0/A_0[t^\infty], (t))\).

    For a Tate ring \(A\), the symbol \(A^\circ\) means the set of all power-bounded elements.
    This gives an open subring of \(A\).
    A subring \(A^+ \subseteq A\) is a \emph{ring of integral elements} if it is open and integrally closed in \(A\) and \(A^+ \subseteq A^\circ\).

\subsection*{Acknowledgment}
I would like to express my gratitude to Kazuma Shimomoto for his time and effort in this research.
Many thanks to Shinnosuke Ishiro for his insightful feedback on the paper's interpretation.
Additionally, we thank Yves Andr\'e, Dimitri Dine, Fumiharu Kato, and Kei Nakazato for reading and providing feedback on this paper.
Finally, I would like to express my deepest gratitude to the anonymous referee.
The benefit of the referee is everywhere.
In particular, the referee gave me the idea of focusing on the ``rigidity lemma'' (\Cref{RigidityLemma}) and precise advice which helps the author to make a clear and readable representation in the Introduction.

\section{Perfectoid Rings} \label{PerfectoidRings}

In this section, we recall and fix some definitions of perfectoid objects.

\begin{definition}[{\cite[Definition 3.5]{bhatt2018Integral}}] \label{PerfectoidRing}
    Let \(S\) be a (non-zero) ring.
    Then \(S\) is an \emph{integral perfectoid ring} if the following conditions hold:
    \begin{enumerate}
        \item There exists an element \(\pi \in S\) such that \(S\) is \(\pi\)-adically complete and \(\pi^p\) divides \(p\) in \(S\).
        \item The Frobenius map \(\map{F}{S/pS}{S/pS}\) is surjective.
        \item The kernel of \(\map{\theta}{\mathbb{A}_{\text{inf}}(S)}{S}\) is principal, where \(\mathbb{A}_{\text{inf}}(S) \defeq W(S^\flat)\).
    \end{enumerate}
    This \(\pi \in S\) is called a \emph{perfectoid element} in this paper.
    Here, we do not require that \(\pi\) is a non-zero-divisor in \(S\).

    Recently, integral perfectoid rings are simply called \emph{perfectoid rings}.
    To avoid confusion with perfectoid Tate rings defined later in \Cref{PerfectoidTateRing}, we do not use the term perfectoid rings, but only integral perfectoid rings and perfectoid Tate rings.
\end{definition}

\begin{lemma}[{(\cite[Lemma 3.9]{bhatt2018Integral})}] \label{pPowerRoots}
    Let \(S\) be an integral perfectoid ring and let \(\pi \in S\) be a perfectoid element.
    Then \(S\) has compatible sequences of \(p\)-power roots \(\{(u\pi)^{1/p^j}\}_{j \geq 0}\) and \(\{(vp)^{1/p^j}\}_{j \geq 0}\) of \(u\pi\) and \(vp\) where \(u\) and \(v\) are unit elements in \(S\).
    
    We fix the element \(\varpi \defeq (vp)^{1/p}\) of \(S\).
    Then \(\varpi\) becomes a perfectoid element of \(S\).
    Without loss of generality, we can assume that a perfectoid element \(\pi\) has a compatible sequence of \(p\)-power roots \(\{\pi^{1/p^j}\}_{j \geq 0}\).
\end{lemma}

\begin{proof}
    The first statement follows from \cite[Lemma 3.9]{bhatt2018Integral}.

    We next check that \(\varpi\) is a perfectoid element of \(S\).
    Note that the conditions (2) and (3) in \Cref{PerfectoidRing} are independent of the choice of a perfectoid element.
    Since \(\varpi^p = vp\) divides \(p\) in \(S\), it suffices to show that \(S\) is \(p\)-adically complete and this is also clear (see \cite[\S 2.1.2]{cesnavicius2023Purity} or \cite[Proposition 2.8]{dospinescu2023Conjecture}).
\end{proof}

\begin{remark}[{\cite[Theorem 3.52]{ishiro2023Perfectoid}}] \label{PerfectoidRingNoWitt}
    Let \(S\) be a ring and let \(\pi\) be an element of \(S\).
    Then \(S\) is an integral perfectoid ring with a perfectoid element \(\pi \in S\) if and only if \(\pi \in S\) satisfies the following:
    \begin{enumerate}
        \item \(S\) is \(\pi\)-adically complete and \(\pi^p\) divides \(p\) in \(S\).
        \item The \(p\)-th power map \(S/\pi S \xrightarrow{a \mapsto a^p} S/\pi^p S\) is an isomorphism of rings.
        \item The multiplicative map
            \begin{align*}
                S[\pi^\infty] & \longrightarrow S[\pi^\infty] \\
                s & \longmapsto s^p
            \end{align*}
            is bijective,
            where the symbol \(S[\pi^\infty]\) is the \emph{ideal of all \(\pi^\infty\)-torsion elements of \(S\)},
            which is defined as
            \begin{equation}
                S[\pi^\infty] = \set{s \in S}{\exists n \in \setZ_{> 0}, \pi^n s = 0 \in S}.
            \end{equation}
    \end{enumerate}
\end{remark}

we recall the definition of semiperfectoid rings.

\begin{definition}[{\cite[Notation 7.1]{bhatt2022Prismsa}}] \label{Semiperfectoid}
    A ring \(S\) is a \emph{semiperfectoid ring} if it is a derived \(p\)-complete ring that is isomorphic to a quotient of an integral perfectoid ring.
\end{definition}

We next explain perfectoid Tate rings. See Notation (\Cref{Notation}) at the end of the Introduction for the basic terminology of Tate rings.

\begin{definition}[\cite{bhatt2018Integral,fontaine2013Perfectoides}] \label{PerfectoidTateRing}
    Let \(A\) be a complete Tate ring (more generally, let \(A\) be a Banach ring).
    Then \(A\) is a \emph{perfectoid Tate ring} if the following conditions hold:
    \begin{enumerate}
        \item \(A\) is uniform, that is, the set of all power-bounded elements \(A^\circ\) is bounded in \(A\).
        \item There exists a pseudo-uniformizer \(\pi \in A\) such that \(\pi^p\) divides \(p\) in \(A^\circ\) and the Frobenius map on \(A^\circ/\pi^p A^\circ\) is surjective.
    \end{enumerate}
    This \(\pi \in A\) is again called a \emph{perfectoid element}.
\end{definition}

The following lemma establishes the connection between integral perfectoid rings and perfectoid Tate rings.

\begin{lemma}[{(\cite[Lemma 3.20]{bhatt2018Integral})}] \label{PerfectoidEquiv}
    Let \(A\) be a Tate ring and let \(A_0^+\) be a ring of integral elements in \(A\).
    If \(A\) is a perfectoid Tate ring, then \(A_0^+\) is an integral perfectoid ring.
    In particular, \(A^\circ\) is an integral perfectoid ring.

    Conversely, if \(A^+_0\) is an integral perfectoid ring and bounded in \(A\), the Tate ring \(A\) is a perfectoid Tate ring.
\end{lemma}

\begin{remark}
    Let \(A\) be a perfectoid Tate ring and let \(\pi \in A\) be a perfectoid element.
    Then \(\pi\) is a non-zero-divisor in \(A\).
    In general, an integral perfectoid ring is not necessarily isomorphic to the set of all power-bounded elements in some perfectoid Tate ring.
\end{remark}

\begin{lemma}[{(\cite[Lemma 3.21]{bhatt2018Integral})}] \label{PerfectoidRingTateRing}
    Let \(S\) be an integral perfectoid ring and let \(\pi \in S\) be a perfectoid element.
    Assume that \(\pi\) is a non-zero-divisor in \(S\).
    Then the Tate ring \(S[1/\pi]\) defined by the pair \((S, (\pi))\) is a perfectoid Tate ring.
    Additionally, \(S \subseteq (S[1/\pi])^{\circ}\), and its cokernel is annihilated by any fractional power of \(\pi\).
\end{lemma}


For convenience, we summarize a brief overview of the properties of perfectoidization.
See \cite{bhatt2022Prismsa,cai2022Perfectoid} for more details.

Let \(S\) be an integral perfectoid ring and let \(R\) be a derived \(p\)-complete \(S\)-algebra.
The \emph{perfectoidization \(R_{\perfd}\) of \(R\)} is defined by using the prismatic cohomology \(\prism_{R/S}\).
Note that \(R_{\perfd}\) is typically a commutative algebra object in \(D^{\geq 0}(S)\) and has a map \(R \to R_{\perfd}\) in \(D(S)\).
The complex \(R_{\perfd}\) is concentrated in degree \(0\) in the following cases:

\begin{itemize}
    \item If \(\Char(S) = p > 0\), \(R_{\perfd}\) coincides with the usual \emph{perfect closure} \(R_{\perf}\) of \(R\).
    \item If \(R\) can be written as a quotient of \(S\), that is, \(R\) is a semiperfectoid ring, then \(R_{\perfd}\) is an integral perfectoid ring and furthermore the map \(R \to R_{\perfd}\) is surjective.
   \item If \(S \to R\) is an integral map, \(R_{\perfd}\) is an integral perfectoid ring.
\end{itemize}

Furthermore, the following property plays an essential role in this paper.

\begin{theorem}[{(\cite[Corollary 8.14]{bhatt2022Prismsa})}] \label{PerfectoidizationUniversal}
    If \(R_{\perfd}\) is concentrated in degree \(0\), it becomes an integral perfectoid ring.
    In this case, the map \(R \to R_{\perfd}\) is the universal map to integral perfectoid rings.
    Namely, every map \(R \to R'\) to an integral perfectoid ring \(R'\) uniquely factors through \(R \to R_{\perfd}\).
\end{theorem}



\section{Uniform Completion} \label{UniformCompletion}

We use the notion of the \emph{uniform completion} of Tate rings.
In this section, we review the uniform completion outlined in \cite{ishizuka2023Mixed}.

\begin{definition}
    A Tate ring \(A\) is \emph{uniform} if the set of all power-bounded elements \(A^\circ\) is bounded.
\end{definition}

Any Tate ring has the structure of a seminormed ring as follows (see \cite[Definition 2.26]{nakazato2022Finite} for more details).

\begin{definition} \label{Seminorm}
Let \(A \defeq A_0[1/t]\) be a Tate ring.
Fix a real number \(c > 1\).
Then, we can define a seminorm \(\map{\norm{\cdot}_{A_0, t, c}}{A}{\setR_{\geq 0}}\) by
\begin{equation*}
    \norm{f}_{A_0, t, c} \defeq \inf_{m \in \setZ} \set{c^m}{t^m f \in A_0}.
\end{equation*}
The seminorm defines a seminormed ring \((A, \norm{\cdot}_{A_0, t, c})\).
The topology of the seminormed ring \((A, \norm{\cdot}_{A_0, t, c})\) is equal to the topology of the Tate ring \(A = A_0[1/t]\).
In particular, the topology induced from the norm \(\norm{\cdot}_{A_0, t, c}\) does not depend on the choices of \(A_0\), \(t\), and \(c\).
So we write the seminorm \(\norm{\cdot}_{A_0, t, c}\) as \(\norm{\cdot}\) for simplicity.

The \emph{spectral seminorm} attached to this seminorm \(\norm{\cdot}\) is defined as
\begin{equation*}
    \norm{f}_{\mathrm{sp}}:=\lim_{n \to \infty} \norm{f^n}^{1/n}.
\end{equation*}
By Fekete's subadditivity lemma, we have \(\norm{f}_{\mathrm{sp}} \leq \norm{f}\).
\end{definition}


\begin{lemma}[{(\cite[Lemma 2.29]{nakazato2022Finite})}] \label{UniformPropeties}
    Let \(A = A_0[1/t]\) be a uniform Tate ring.
    Then \(A^\circ\) is equal to the unit disk \(A_{\norm{\cdot}_{\mathrm{sp}} \leq 1}\) by the spectral seminorm in \(A\).
\end{lemma}

Next, we define the uniformization and uniform completion of Tate rings as in \cite[Exercise 7.2.6]{bhatt2017Lecture} and \cite[Definition 2.8.13]{kedlaya2015Relative}.

\begin{definition} \label{UniformizationUniformCompletionWRT}
    Let \(A\) be a Tate ring.
    Fix a pair of definition \((A_0, (t))\) of \(A\) and a ring of integral elements \(A_0^+\) of \(A\) such that \(A_0 \subseteq A_0^+\).
    We define the following terminology.
    \begin{enumerate}
        \item The \emph{uniformization of \(A\) with respect to \((A_0, (t))\) and \(A_0^+\)} is the Tate ring \(A_0^+[1/t]\).
        \item The \emph{uniform completion of \(A\) with respect to \((A_0, (t))\) and \(A_0^+\)} is the completion of the Tate ring \(A_0^+[1/t]\). This is in fact the Tate ring \(\widehat{A_0^+}[1/t]\) where \(\widehat{A_0^+}\) is the \(t\)-adic completion of \(A_0^+\).
    \end{enumerate}
\end{definition}

\begin{remark} \label{UniformizationUniformCompletionIndep}
    At first glance, the above definitions depend on the choices of \((A_0, (t))\) and \(A_0^+\). The motivation for these definitions is that we wanted to define them ``functorially''.

    Furthermore, even if we take a different pair of definition \((A_0', (t'))\) of \(A\) and a different ring of integral elements \(A_0^{'+}\) such that \(A_0 \subseteq A_0^{'+}\), the uniformization of \(A\) with respect to \((A_0', (t'))\) and \(A_0^{'+}\) is isomorphic to the uniformization of \(A\) with respect to \((A_0, (t))\) and \(A_0^+\).
    Its isomorphism is obtained by the identity map on the abstract ring \(A = A_0^{'+}[1/t] = A_0^+[1/t]\) (see \cite[Lemma 5.5]{ishizuka2023Mixed} and \cite[Lemma 2.3]{nakazato2022Finite}).
    In particular, the same statement is true for the uniform completion.
    So the next definitions are well-defined.

\end{remark}

\begin{definition} \label{UniformizationUniformCompletion}
    Let \(A\) be a Tate ring.
    The \emph{uniformization \(A^u\) of \(A\)} (resp., \emph{uniform completion \(A^{\widehat{u}}\) of \(A\)}) is the uniformization (resp., uniform completion) of \(A\) with respect to a pair of definition \((A_0, (t))\) and a ring of integral elements \(A_0^+\) of \(A\) such that \(A_0 \subseteq A_0^+\).
    By the above \Cref{UniformizationUniformCompletionIndep}, these definitions are independent of the choices of \((A_0, (t))\) and \(A_0^+\).

\end{definition}


Recall that the canonical map of Tate rings \(i \colon A \to A^u \to A^{\widehat{u}}\) has the following universal property.

\begin{proposition}[{(\cite[Proposition 5.6]{ishizuka2023Mixed})}] \label{UniversalityUniformCompletion}
    Let \(A\) be a Tate ring.
    Then the uniform completion \(A^{\widehat{u}}\) is a uniform complete Tate ring.
    Furthermore, the canonical map \(\map{i}{A}{A^{\widehat{u}}}\) is the universal map to uniform complete Tate rings.
    That is, every map of Tate rings \(\map{h}{A}{B}\), where \(B\) is a uniform complete Tate ring uniquly factors through \(\map{i}{A}{A^{\widehat{u}}}\).
\end{proposition}


    

We record some lemmas as follows.

\begin{lemma} \label{PowerBoundedUniformization}
    Let \(A\) be a Tate ring.
    Then, the completion \(\widehat{A^{u \circ}}\) of the set of all power-bounded elements \(A^{u \circ}\) of \(A^u\) is isomorphic to the set of all power-bounded elements \((A^{\widehat{u}})^\circ\) of \(A^{\widehat{u}}\) as topological ring.
\end{lemma}

\begin{proof}
    Fix a pair of definition \((A_0, (t))\) of \(A\) and a ring of integral elements \(A_0^+\) of \(A\) such that \(A_0 \subseteq A_0^+\).
    By \cite[Lemma 2.3]{nakazato2022Finite}, we have an inclusion \(t (A^+_0)^*_A \subseteq A_0^+\) where \((A_0^+)^*_A\) is the complete integral closure of \(A_0^+\) in \(A\).
    By \cite[Proposition 2.4]{nakazato2022Finite}, the canonical map \(A^{\widehat{u}} = \widehat{A^+_0}[1/t] \to \widehat{(A_0^+)^*_{A^u}}[1/t]\) is an isomorphism of Tate rings and this map induces an isomorphism of topological rings
    \(\parenlr{\widehat{A^+_0}}^*_{A^{\widehat{u}}} \to \widehat{(A_0^+)^*_{A^u}}\).
    Since \(A^{\widehat{u}}\) (resp, \(A^u\)) has a ring of definition \(\widehat{A^+_0}\) (resp, \(A_0^+\)) and the complete integral closure is equal to the set of all power-bounded elements by \cite[Lemma 2.13]{nakazato2022Finite}, we have an isomorphism of topological rings \((A^{\widehat{u}})^\circ \xrightarrow{\cong} \widehat{A^{u \circ}}\).
\end{proof}

\begin{lemma}[{(\cite[Proposition 5.6]{ishizuka2023Mixed})}] \label{UniformCompletionAlmostElement}
    Assume that a Tate ring \(A = A_0[1/t]\) has a compatible sequence of \(p\)-power roots \(\{t^{1/p^j}\}_{j \geq 0}\) of \(t\).
    The inclusion map \((A^{\widehat{u}})^\circ \hookrightarrow ((A^{\widehat{u}})^\circ)_* \defeq t^{-1/p^\infty} ((A^{\widehat{u}})^\circ)\) is an isomorphism of rings.
\end{lemma}

\begin{proof}
    For any \(f \in ((A^{\widehat{u}})^\circ)_*\), we have \(t^{1/p^n} f \in (A^{\widehat{u}})^\circ\) and then \(t f^{p^n} \in (A^{\widehat{u}})^\circ\) for any \(n \in \setZ_{> 0}\).
    Since \((A^{\widehat{u}})^\circ\) is a ring of definition of \(A^{\widehat{u}}\), we have \(\norm{f^{p^n}} \leq c\) for a fixed \(c > 1\) by \Cref{Seminorm}.
    In particular, \(\norm{f^{p^n}}^{1/p^n} \leq c^{1/p^n}\) for any \(n \in \setZ_{> 0}\).
    Taking the limit \(n \to \infty\), we have \(\norm{f}_{\mathrm{sp}} \leq \lim_{n \to \infty} c^{1/p^n} = 1\).
    This shows the inclusion \(((A^{\widehat{u}})^\circ)_* \subseteq (A^{\widehat{u}})_{\norm{\cdot}_{\mathrm{sp}} \leq 1} = (A^{\widehat{u}})^\circ\) by \Cref{UniformPropeties}.
\end{proof}

\section{Some Ring-Theoretic Properties of Pre-Perfectoid Pairs} \label{PrePerfectoidPair}

Let \(R\) be an integral perfectoid ring and let \(\pi \in R\) be a perfectoid element.
Our goal in this section is to convert a situation where \(\pi\) is a zero-divisor into a situation where it is a non-zero-divisor, following the approach of \cite[\S 2.3.2]{andre2020Weak}.
Our argument is based on \cite{andre2018Lemme,andre2020Weak} and is similar to \cite{scholze2012Perfectoida,bhatt2017Lecture}.
If it is sufficient to consider only \(\pi\)-torsion-free rings (resp., \(p\)-torsion-free rings), the symbol \((-)^{\pi \mathrm{tf}}\) defined in \Cref{DefPrePerfectoidPair} (resp., \((-)^{\ptf}\)) can be removed.

For the sake of generality, we define pre-perfectoid pairs as follows.

\begin{definition} \label{DefPrePerfectoidPair}
    Let \((S, \pi)\) be a pair such that \(S\) is a ring and \(\pi\) is an element of \(S\) which has a compatible sequence of \(p\)-power roots \(\{\pi^{1/p^j}\}_{j \geq 0}\) in \(S\) and \(\pi^p\) divides \(p\) in \(S\).
    If the \(p\)-th power map \(S/\pi S \xrightarrow{a \mapsto a^p} S/\pi^p S\) is isomorphism, we call such a pair \((S, \pi)\) \emph{pre-perfectoid pair}.

    For a pre-perfectoid pair \((S, \pi)\), the \emph{\(\pi\)-torsion-free quotient \(S^{\pi \mathrm{tf}}\) of \(S\)} is defined as the quotient ring
    \begin{equation*}
        S^{\pi \mathrm{tf}} \defeq S/S[\pi^\infty] \cong \Image(S \to S[1/\pi]),
    \end{equation*}
    where \(S[\pi^\infty]\) is the ideal of all \(\pi^\infty\)-torsion elements of \(S\).
    Note that \(S \twoheadrightarrow S^{\pi \mathrm{tf}}\) is an isomorphism if and only if \(\pi\) is a non-zero-divisor of \(S\).
    In the case of \(S[1/\pi] = S[1/p]\), \(S^{\pi \mathrm{tf}}\) is equal to the \(p\)-torsion-free quotient \(S^{\ptf}\) of \(S\) defined in \Cref{pTorsionCase}.
\end{definition}

For example, an integral perfectoid ring \(R\) and a perfectoid element \(\pi\) of \(R\) form a pre-perfectoid pair \((R, \pi)\) because of \Cref{PerfectoidRingNoWitt} (2).

\begin{definition}
    Let \((S, \pi)\) be a pre-perfectoid pair.
    An \(S\)-module \(M\) is called \emph{\((\pi)^{1/p^\infty}\)-almost zero} if \(\pi^{1/p^n} \cdot M = 0\) for any \(n \in \setZ_{> 0}\).
    A map of \(S\)-modules \(N \to M\) is called \emph{\((\pi)^{1/p^\infty}\)-almost injective (resp, surjective)} if its kernel (resp, cokernel) is \((\pi)^{1/p^\infty}\)-almost zero.
    If a map of \(S\)-modules \(N \to M\) is \((\pi)^{1/p^\infty}\)-almost injective and surjective, we call the map a \emph{\((\pi)^{1/p^\infty}\)-almost isomorphism}.
\end{definition}

\begin{definition} \label{AlmostElement}
    Let \((S, \pi)\) be a pre-perfectoid pair.
    We define  \emph{the set of almost elements of \(S\)} (see for example \cite[Lemma 5.3]{scholze2012Perfectoida}):
    \begin{equation}
        S_* \defeq \pi^{-1/p^\infty} S \defeq \pi^{-1/p^\infty} S^{\pi \mathrm{tf}} \defeq \set{s \in S[1/\pi]}{\forall n \in \setZ_{> 0}, \pi^{1/p^n} s \in S^{\pi \mathrm{tf}} \subseteq S[1/\pi]}.
    \end{equation}

    For any \(S\)-module \(M\), \emph{the set of almost elements of \(M\)} is defined as
    \begin{equation}
        M_* \defeq \pi^{1/p^\infty} M \defeq \set{m \in M[1/\pi]}{\forall n \in \setZ_{> 0}, \pi^{1/p^n} m \in M^{\pi \mathrm{tf}} \subseteq M[1/\pi]}
    \end{equation}
    where \(M^{\pi \mathrm{tf}}\) is the \(\pi\)-torsion-free quotient of \(M\), that is, the image of \(M\) in \(M[1/\pi]\).
\end{definition}

\begin{remark} \label{TopRingStr}
    Let \((S, \pi)\) be a pre-perfectoid pair.
    Note that \(S[1/\pi]\) and \(S_*[1/\pi]\) are isomorphic as abstract rings.
    Furthermore, as Tate rings, \(S[1/\pi]\) has a pair of definition \((S^{\pi \mathrm{tf}}, (\pi))\)
    and \(S_*[1/\pi]\) has a pair of definition \((S_*, (\pi))\).
    Because of \(\pi S_* \subseteq S^{\pi \mathrm{tf}} \subseteq S_*\), these are isomorphic as topological rings.
\end{remark}

\begin{lemma} \label{AlmostElementAlmostIsom}
    Let \((S, \pi)\) be a pre-perfectoid pair.
    Assume that \(S\) is reduced (e.g. \(S\) is an integral perfectoid ring, see \cite[\S 2.1.3]{cesnavicius2023Purity} or \cite[Proposition 2.20]{dospinescu2023Conjecture}).
    Then, the canonical map \(S \to S_*\) is a \((\pi)^{1/p^\infty}\)-almost isomorphism.
\end{lemma}

\begin{proof}
    The kernel of \(S \to S[1/\pi]\) is isomorphic to the ideal \(S[\pi^\infty]\) of all \(\pi^\infty\)-torsion elements of \(S\).
    Since \(S\) is reduced, \(S[\pi^\infty]\) is \((\pi)^{1/p^\infty}\)-almost zero.
    In particular, \(S \to S^{\pi \mathrm{tf}}\) is (usual) surjective and \((\pi)^{1/p^\infty}\)-almost injective.
    Furthermore, the inclusion \(S^{\pi \mathrm{tf}} \subseteq \pi^{-1/p^\infty} S\) in \(S[1/\pi]\) is (usual) injective and \((\pi)^{1/p^\infty}\)-almost surjective by definition.
    This completes the proof.
\end{proof}

\begin{lemma} \label{FrobAlmostIsom}
    Let \((S, \pi)\) a pre-perfectoid pair.
    Then the \(p\)-th power map \(S_*/\pi S_* \xrightarrow{a \mapsto a^p} S_*/\pi^p S_*\) is injective and \((\pi)^{1/p^\infty}\)-almost surjective (later in \Cref{AlmostFrobSurj}, this will become a (usual) surjective map).
    In particular, \(S_*\) is \(p\)-root closed in \(S_*[1/\pi]\). Namely, if \(x \in S_*[1/\pi]\) satisfies \(x^{p^n} \in S_*\) for some \(n \in \setZ_{> 0}\), then \(x \in S_*\).
\end{lemma}

\begin{proof}
    If \(\pi\) is a non-zero-divisor of \(S\), this lemma is similar to \cite[Proposition 2.16 (b)]{dospinescu2023Conjecture}.
    We only have to make the same proof as \cite[Lemma 5.6]{scholze2012Perfectoida}, being careful that \(S\) is not necessarily \(\pi\)-torsion-free in our case.
    
    First, we show that \(S_*/\pi S_* \xrightarrow{a \mapsto a^p} S_*/\pi^p S_*\) is injective.
    Let \(t \in S_*\) be an element in the kernel of the \(p\)-th power map.
    There exists some \(t' \in S_*\) such that \(t^p = \pi^p t'\) in \(S_* \subseteq S[1/\pi]\).
    By definition of \(S_*\), multiplying \(\pi^{1/p^n}\) by the equation for each \(n \in \setZ_{> 0}\), we have
    \begin{equation} \label{eq1}
        (\pi^{1/p^{n+1}} t)^p = \pi^{1/p^n} t^p = \pi^p (\pi^{1/p^n} t') \in \pi^p S^{\pi \mathrm{tf}}.
    \end{equation}
    Moreover, \(\pi^{1/p^{n+1}} t\) and \(\pi^{1/p^n} t'\) are elements of \(S^{\pi \mathrm{tf}}\).
    Then, there exist some elements \(s_{n+1}\) and \(s'_n\) in \(S\) such that \(s_{n+1}/1 = \pi^{1/p^{n+1}} t\) and \(s'_n/1 = \pi^{1/p^n} t'\) in \(S^{\pi \mathrm{tf}} \subseteq S[1/\pi]\).
    By the above equation (\ref{eq1}), we have \(s_{n+1}^p/1 = \pi^p s'_n/1\) in \(S^{\pi \mathrm{tf}} \subseteq S[1/\pi]\) and thus, \(s_{n+1}^p - \pi^p s'_n\) is in \(S[\pi^{\infty}] \subseteq S\), which is \((\pi)^{1/p^\infty}\)-almost zero.
    For any \(m \in \setZ_{> 0}\), we have
    \begin{equation*}
        S \ni 0 = \pi^{1/p^m} (s_{n+1}^p - \pi^p s'_n)
        = (\pi^{1/p^{m+1}} s_{n+1})^p - \pi^p (\pi^{1/p^m} s'_n).
    \end{equation*}
    In particular, \(\pi^{1/p^{m+1}} s_{n+1}\) is in the kernel of the \(p\)-th power map \(S/\pi S \to S/\pi^p S\), which is zero by assumption of \((S, \pi)\), and so \(\pi^{1/p^{m+1}} s_{n+1}\) is in \(\pi S\).
    Passing to \(S^{\pi \mathrm{tf}} \subseteq S[1/\pi]\), we have
    \begin{equation*}
        \pi S^{\pi \mathrm{tf}} \ni \pi^{1/p^{m+1}} s_{n+1}/1 = \pi^{1/p^{m+1}} \pi^{1/p^{n+1}} t
    \end{equation*}
    for any \(n, m \in \setZ_{> 0}\).
    Then, \(t\) is in \((\pi S^{\pi \mathrm{tf}})_*\), which is defined in \Cref{AlmostElement}.
    The next equality (\ref{CommutesAlmostElement}) shows the injectivity of \(S_*/\pi S_* \xrightarrow{a \mapsto a^p} S_*/\pi^p S_*\):
    \begin{equation} \label{CommutesAlmostElement}
        (\pi S^{\pi \mathrm{tf}})_* = \pi (S^{\pi \mathrm{tf}})_* = \pi S_* \subseteq S[1/\pi].
    \end{equation}
    

    \begin{proof}[Proof of (\ref{CommutesAlmostElement})]
        By definition, we have \((S^{\pi \mathrm{tf}})_* = S_*\) and then the second equality is clear.
        Any element of \(\pi (S^{\pi \mathrm{tf}})_*\) can be written as \(\pi t\) by using some element \(t \in (S^{\pi \mathrm{tf}})_*\).
        Because of \(\pi^{1/p^n} (\pi t) = \pi (\pi^{1/p^n} t) \in \pi S^{\pi \mathrm{tf}}\) for any \(n \in \setZ_{> 0}\), we have \(\pi t \in (\pi S^{\pi \mathrm{tf}})_*\) and so \(\pi (S^{\pi \mathrm{tf}})_* \subseteq (\pi S^{\pi \mathrm{tf}})_*\).

        Conversely, take any element \(x \in (\pi S^{\pi \mathrm{tf}})_* \subseteq S[1/\pi]\).
        Then, there exists \(t_n \in S^{\pi \mathrm{tf}}\) such that \(\pi^{1/p^n} x = \pi t_n \in \pi S^{\pi \mathrm{tf}}\) for each \(n \in \setZ_{> 0}\).
        Since \(S[1/\pi]\) is \(\pi^{1/p^n}\)-torsion-free, we have
        \begin{equation}
            S[1/\pi] \ni \pi^{1/p^n} (x/\pi) = t_n \in S^{\pi \mathrm{tf}}
        \end{equation}
        for any \(n \in \setZ_{> 0}\).
        This shows that \(x/\pi \in S[1/\pi]\) is in \((S^{\pi \mathrm{tf}})_*\) and thus \(x = \pi (x/\pi)\) is in \(\pi (S^{\pi \mathrm{tf}})_*\).
    \end{proof}

    Second, we show that \(S_*/\pi S_* \xrightarrow{a \mapsto a^p} S_*/\pi^p S_*\) is \((\pi)^{1/p^\infty}\)-almost surjective.
    Take any element \(x \in S_* \subseteq S[1/\pi]\).
    For each \(n \in \setZ_{> 0}\), we have \(\pi^{1/p^n} x \in S^{\pi \mathrm{tf}}\) and thus there exists an \(s_n \in S\) such that \(\pi^{1/p^n} x = s_n/1 \in S^{\pi \mathrm{tf}}\).
    Since \(S/\pi S \xrightarrow{a \mapsto a^p} S/\pi^p S\) is surjective, there exists some \(s'_n \in S\) such that \((s'_n)^p - s_n\) is in \(\pi^p S\).
    Then, we have
    \begin{equation}
        \pi^{1/p^n} x - (s'_n/1)^p = s_n/1 - (s'_n)^p/1 \in \pi^p S^{\pi \mathrm{tf}} \subseteq \pi^p S_*
    \end{equation}
    and \(s'_n/1\) is in \(S^{\pi \mathrm{tf}} \subseteq S_*\).
    This shows that \(S_*/\pi S_* \xrightarrow{a \mapsto a^p} S_*/\pi^p S_*\) is \((\pi)^{1/p^\infty}\)-almost surjective.

    Finally, since \(S_*\) is \(\pi\)-torsion-free and \(S_*/\pi S_* \xrightarrow{a \mapsto a^p} S_*/\pi^p S_*\) is injective as above, by \cite[(2.1.7.1)]{cesnavicius2023Purity}, we can show that \(S_*\) is \(p\)-root closed in \(S_*[1/\pi]\).
\end{proof}

\begin{lemma} \label{UniformalityAlmostElement}
    Let \((S, \pi)\) be a pre-perfectoid pair.
    Then \(S_*[1/\pi]\) is a uniform Tate ring which satisfies \((S_*[1/\pi])^\circ = S_*\).
    In particular, \(S_*\) is isomorphic to the unit disk of \(S_*[1/\pi]\) with respect to the spectral seminorm induced from the seminorm of \(S_*[1/\pi]\) as defined in \Cref{Seminorm}.
\end{lemma}

\begin{proof}
    
    Proceeding as in \cite[Lemma 5.6]{scholze2012Perfectoida}, the \(p\)-root closedness of \(S_*\) in \(S_*[1/\pi]\) proved in \Cref{FrobAlmostIsom} shows that \((S_*[1/\pi])^\circ = S_*\).
    In fact, any element \(x \in (S_*[1/\pi])^\circ\) makes a topologically nilpotent element \(\pi^{1/p^n} x\) for each \(n \in \setZ_{> 0}\).
    Thus there exists \(N = N(n) \in \setZ_{> 0}\) such that \((\pi^{1/p^n} x)^{p^N} \in S_*\).
    This shows the equality above and thus \(S_*[1/\pi]\) is a uniform Tate ring.
    Moreover, \Cref{UniformPropeties} shows that
    \begin{equation}
        S_* = (S_*[1/\pi])^\circ = (S_*[1/\pi])_{\norm{\cdot}_{\mathrm{sp}} \leq 1}
    \end{equation}
    and we finish the proof.
\end{proof}

\begin{corollary} \label{AlmostFrobSurj}
    Let \((S, \pi)\) be a pre-perfectoid pair.
    Stronger than \Cref{FrobAlmostIsom}, we can show that the \(p\)-th power map \(S_*/\pi S_* \xrightarrow{a \mapsto a^p} S_*/\pi^p S_*\) is (usual) surjective.
\end{corollary}

\begin{proof}
    This proof is similar to \cite[Lemma 4.2]{ishizuka2023Mixed}.
    Fix an element \(y \in S_*\).
    Since the \(p\)-th power map is \((\pi)^{1/p^\infty}\)-almost surjective by \Cref{FrobAlmostIsom}, there exist elements \(a\) and \(b\) in \(S_*\) such that \(\pi y = a^p + \pi^p b \in S_*\).
    Set \(z \defeq a/\pi^{1/p} \in S_*[1/\pi]\). This satisfies
    \begin{equation}
        z^p = a^p/\pi = y - \pi^{p-1} b \in S_*.
    \end{equation}
    Since \(S_*\) is \(p\)-root closed in \(S_*[1/\pi]\) by \Cref{FrobAlmostIsom}, we can show that
    \(z \in S_*\).
    The equality \(\pi y = \pi z^p + \pi^p b \in S_*\) and \(\pi\)-torsion-freeness of \(S_* \subseteq S_*[1/\pi]\) show that \(y = z^p + \pi^{p-1} b \in S_*\).
    In particular, every element of \(S_*/\pi^{p-1} S_*\) and \(S_*/\pi S_*\) is a \(p\)-th power.
    The inclusion \(\pi S_* \subseteq \pi^{(p-1)/p} S_*\) induces the commutative diagram
    \begin{center}
        \begin{tikzcd}
            0 \arrow[r] & S_*/\pi^{1/p} S_* \arrow[r, "\cdot \pi^{(p-1)/p}"] \arrow[d, two heads] & S_*/\pi S_* \arrow[r] \arrow[d] & S_*/\pi^{(p-1)/p} S_* \arrow[r] \arrow[d, two heads] & 0 \\
            0 \arrow[r] & S_*/\pi S_* \arrow[r, "\cdot \pi^{p-1}"]                                & S_*/\pi^p S_* \arrow[r]         & S_*/\pi^{p-1} S_* \arrow[r]                          & 0
        \end{tikzcd}
    \end{center}
    where the rows are exact sequences and the vertical maps are the \(p\)-th power maps.
    Since the right-most map and left-most map are surjective, the middle map is also surjective by the five lemma.
\end{proof}

Finally, we show the stability of pre-perfectoid pairs under taking \((-)_*\)
and the generalization of \Cref{PerfectoidRingTateRing}.

\begin{theorem} \label{UniformPerfectoid}
    For any pre-perfectoid pair \((S, \pi)\), the induced pair \((S_*, \pi)\) is again a pre-perfectoid pair.

    In particular, for any integral perfectoid ring \(R\) and a perfectoid element \(\pi \in R\) (not necessarily a non-zero-divisor),
    \(R_* = \pi^{-1/p^\infty} R\) is a \(\pi\)-torsion-free integral perfectoid ring with a perfectoid element \(\pi\).

\end{theorem}

\begin{proof}
    By \Cref{FrobAlmostIsom} and \Cref{AlmostFrobSurj}, the pair \((S_*, \pi)\) is a pre-perfectoid pair.

    Assume that \(R\) is an integral perfectoid ring.
    By the above paragraph of \cite[(2.1.3.1)]{cesnavicius2023Purity}, the image \(R^{\pi \mathrm{tf}}\) of \(R\) in \(R[1/\pi]\) is a \(\pi\)-torsion-free integral perfectoid ring.
    In particular, similar to the proof of \cite[Lemma 5.6]{scholze2012Perfectoida}, \(R_* = (R^{\pi \mathrm{tf}})_*\) is also an integral perfectoid ring with a perfectoid element \(\pi\).
\end{proof}

\begin{remark} \label{PerfectoidPair}
    By \Cref{pPowerRoots}, any integral perfectoid ring \(R\) has a compatible sequence of \(p\)-power roots \(\{\varpi^{1/p^j}\}_{j \geq 0}\) of \(\varpi \in R\) such that \(\varpi^p\) is some unit multiple of \(p\) in \(R\).
    Then \(R\) forms a pre-perfectoid pair \((R, \varpi)\).
    In particular, this pre-perfectoid pair \((R, \varpi)\) satisfies all the statements in this section.
\end{remark}

\section{Calculations of the Perfectoidization}

We next show the ``universality'' of uniformizations and uniform completions for pre-perfectoid pairs and deduce our main theorem (\Cref{Maintheorem2}).
As in the previous section, we can get rid of the symbol \((-)^{\ptf}\) if we only consider \(p\)-torsion-free rings.

\begin{proposition} \label{UniversalityUniform}
    Let \(A_0\) be a \(\pi\)-adically topological ring for an element \(\pi \in A_0\).
    Set a Tate ring \(A \defeq A_0[1/\pi]\).
    Let \(S\) be an \(A_0\)-algebra such that the image of \(\pi\) in \(S\) is a non-zero element.
    Assume that \(S\) is an integral perfectoid ring and \(\pi\) is a perfectoid element of \(S\).

    Then, there exists a unique map of topological rings \((A^{\widehat{u}})^\circ \to S_*\) such that the following diagram commutes:
    \begin{center}
        \begin{tikzcd}
            A_0 \arrow[r] \arrow[d]   & S \arrow[d] \\
            (A^{\widehat{u}})^\circ \arrow[r, "\exists !", dashed] & S_*.
        \end{tikzcd}
    \end{center}
\end{proposition}

\begin{proof}
    Since the \(A_0\)-algebra structure \(A_0 \to S\) is a continuous map with respect to the \(\pi\)-adic topology, this can be extended to the map of Tate rings \(A \to S[1/\pi]\).
    By \Cref{TopRingStr} and \Cref{UniformalityAlmostElement}, the Tate ring \(S[1/\pi]\) is isomorphic to the complete uniform Tate ring \(S_*[1/\pi]\).
    
    The universality of uniform completion (\Cref{UniversalityUniformCompletion}) gives a unique map of Tate rings \(A^{\widehat{u}} \to S_*[1/\pi]\) which extends \(A \to S[1/\pi] = S_*[1/\pi]\).
    Taking the set of all power-bounded elements, we have a map of topological rings \((A^{\widehat{u}})^\circ \to S_*\) such that the above diagram commutes
    by \Cref{UniformalityAlmostElement}.
    The uniqueness of this map is clear.
\end{proof}

To prove our main theorem (\Cref{Maintheorem2}), we state the next lemma which is in the form of ``rigidity lemma'' as in \Cref{RigidityLemma}.

\begin{lemma} \label{IdentityMap}
    Let \(A_0\) be a derived \(p\)-complete algebra over some integral perfectoid ring of characteristic \(0\).\footnote{Note that any integral perfectoid ring does not contain \(\setQ\). So a (\(p\)-adically complete) integral perfectoid ring \(R\) is of characteristic \(0\) means that \(R\) contains \(\setZ\) as a subring and \(pR \neq R\). We do not assume that \(R\) is \(p\)-torsion-free.}
    Assume that \((A_0)_{\perfd}\) is an (honest) integral perfectoid ring.
    Then any \(A_0\)-algebra map \((A_0)_{\perfd, *} \to (A_0)_{\perfd, *}\)\footnote{Recall that an integral perfectoid ring \((A_0)_{\perfd}\) has a compatible sequence \(\{\varpi^{1/p^n}\}_{n \geq 0}\) of \(p\)-power roots of \(\varpi\) such that \(\varpi^p\) is a unit multiple of \(p\) in \((A_0)_{\perfd}\) (see \Cref{pPowerRoots} and \Cref{PerfectoidPair}). So \((A_0)_{\perfd, *}\) is the set of almost elements of the pre-perfectoid pair \(((A_0)_{\perfd}, \varpi)\). That is, \((A_0)_{\perfd, *} = \varpi^{-1/p^\infty} (A_0)_{\perfd} \subseteq (A_0)_{\perfd}[1/\varpi] = (A_0)_{\perfd}[1/p]\).} is the identity map.
\end{lemma}

\begin{proof}
    Take any map \(\map{f}{(A_0)_{\perfd,*}}{(A_0)_{\perfd,*}}\) such that the following diagram commutes
    \begin{center}
        \begin{tikzcd}
            (A_0)_{\perfd} \arrow[d, "c"]      & A_0 \arrow[r] \arrow[l] & (A_0)_{\perfd} \arrow[d, "c"] \\
            {(A_0)_{\perfd,*}} \arrow[rr, "f"] &                         & {(A_0)_{\perfd,*}}           
        \end{tikzcd}
    \end{center}
    where the vertical maps are canonical ones.
    Since \((A_0)_{\perfd,*}\) is also an integral perfectoid ring by \Cref{UniformPerfectoid} and \((A_0)_{\perfd}\) is a universal integral perfectoid ring over \(A_0\), the composite \(f \circ c\) is nothing but the canonical map \(\map{c}{(A_0)_{\perfd}}{(A_0)_{\perfd,*}}\).
    For any \(x \in (A_0)_{\perfd,*}\), there exists an element \(a \in (A_0)_{\perfd}\) such that \(c(a) = p x\) in \((A_0)_{\perfd,*}\).
    Then, we have \(p f(x) = f(p x) = f(c(a)) = c(a) = p x\) in \((A_0)_{\perfd,*}\).
    Since \((A_0)_{\perfd,*}\) is \(p\)-torsion-free, we have \(f(x) = x\) and we are done.
\end{proof}

As a consequence of this lemma, we can show the following relation between \((-)^{\ptf}\) and \((-)_{\perfd}\).

\begin{corollary} \label{Interchangeable}
    Let \(A_0\) be a derived \(p\)-complete algebra over some integral perfectoid ring.
    Assume that \((A_0)_{\perfd}\) and \((A_0^{\ptf})_{\perfd}\) are (honest) integral perfectoid rings.
    Then the canonical map \((A_0)_{\perfd} \to (A_0^{\ptf})_{\perfd}\) induces the isomorphism
    \begin{equation*}
        ((A_0)_{\perfd})^{\ptf} \xrightarrow{\cong} (A_0^{\ptf})_{\perfd}.
    \end{equation*}
\end{corollary}

    Since the above isomorphism prevents confusion when \(((A_0)_{\perfd})^{\ptf}\) is written as \((A_0)_{\perfd}^{\ptf}\), we use this symbol \((-)_{\perfd}^{\ptf}\) in the following.

\begin{proof}[Proof of \Cref{Interchangeable}]
    The canonical surjective map \(A_0 \twoheadrightarrow A_0^{\ptf}\) induces a map \((A_0)_{\perfd} \to (A_0^{\ptf})_{\perfd}\).
    Since \((A_0^{\ptf})_{\perfd}\) is \(p\)-torsion-free by \cite[Lemma A.2]{ma2022Analogue}, we have a unique map \(\map{\varphi}{((A_0)_{\perfd})^{\ptf}}{(A_0^{\ptf})_{\perfd}}\) which extends the map \((A_0)_{\perfd} \to (A_0^{\ptf})_{\perfd}\).
    
    Conversely, by \cite[Proposition 2.19]{dospinescu2023Conjecture}, the \(p\)-torsion-free quotient \(((A_0)_{\perfd})^{\ptf}\) is also an integral perfectoid ring over \(A_0\) and thus there exists a unique map \(A_0^{\ptf} \to ((A_0)_{\perfd})^{\ptf}\) which extends the structure map \(A_0 \to ((A_0)_{\perfd})^{\ptf}\).
    The universality of perfectoidization (\Cref{PerfectoidizationUniversal}) shows that this map can be extended to a map \(\map{\psi}{(A_0^{\ptf})_{\perfd}}{((A_0)_{\perfd})^{\ptf}}\) uniquely.
    
    Then the \(A_0^{\ptf}\)-algebra map \(\varphi \circ \psi\) is the identity map because of the universality of perfectoidizations (\Cref{PerfectoidizationUniversal}).
    On the other hand, the \((A_0)_{\perfd}\)-algebra map \(\psi \circ \varphi\) can be extended to the \((A_0)_{\perfd}\)-algebra endomorphism on \((((A_0)_{\perfd})^{\ptf})_*\) because of the inclusion \(((A_0)_{\perfd})^{\ptf} \subseteq (((A_0)_{\perfd})^{\ptf})_* \subseteq ((A_0)_{\perfd})[1/p]\).
    Then the above rigidity lemma (\Cref{IdentityMap}) shows that \(\psi \circ \varphi\) is in fact the identity map on \(((A_0)_{\perfd})^{\ptf}\).
    This completes the proof.
\end{proof}

In the proof of \cite[Theorem 4.4]{dine2022Topologicala}, Dine states that a quotient of a perfectoid Tate ring by some ideal has the perfectoidization that is isomorphic to its uniform completion.
We reformulate the proof for the situation of integral perfectoid rings as follows.

\begin{theorem}[(cf. \cite{dine2022Topologicala})] \label{MainTheorem1}
    Let \(A_0\) be a derived \(p\)-complete algebra over some integral perfectoid ring of characteristic \(0\).
    Set a Tate ring \(A \defeq A_0[1/p]\).
    Assume that the perfectoidization \((A_0)_{\perfd}\) is an (honest) integral perfectoid ring and the uniform completion \(A^{\widehat{u}}\) of \(A\) is a perfectoid Tate ring.
    
    Then \(A^{\widehat{u}}\) is isomorphic to \((A_0)_{\perfd}[1/p]\) as a Tate ring.
    In particular, \((A^{\widehat{u}})^\circ\) and \((A_0)_{\perfd, *}\) are isomorphic as rings.

\end{theorem}

\begin{proof}
    If we know that \(A^{\widehat{u}}\) is a perfectoid Tate ring, the same proof of \cite[Theorem 4.4]{dine2022Topologicala} induces the isomorphism \(A^{\widehat{u}} \cong (A_0)_{\perfd}[1/p]\) by checking the universality of perfectoidizations and uniform completions as in \Cref{UniversalityUniform}.
    By \Cref{UniformalityAlmostElement}, taking the set of all power-bounded elements induces the isomorphism \((A^{\widehat{u}})^\circ \cong (A_0)_{\perfd, *}\).
\end{proof}

\begin{remark} \label{UniformCompletionPerfectoid}
    Note that \((A^{\widehat{u}})^\circ\) is an integral perfectoid ring if and only if \(A^{\widehat{u}}\) is a perfectoid Tate ring by \Cref{PerfectoidEquiv}.
    So the conditions for \(A^{\widehat{u}}\) to be a perfectoid Tate ring are studied in \cite[Theorem 3.3.18 (ii)]{kedlaya2019Relative}.
    For example, if \(A_0\) is a semiperfectoid ring, then \((A^{\widehat{u}})^\circ\) is an integral perfectoid ring as shown in \cite[Theorem 4.4]{dine2022Topologicala}.
    Compare \Cref{pRootClosurePerfectoid} below.
\end{remark}

We recall that taking the set of all power-bounded elements is not only a topological operation but also an algebraic operation as mentioned in the next lemma.
This lemma is used in the proof of \Cref{ExplicitPerfectoidization}.

\begin{lemma} \label{UniformCompRealization}
    Let \(A_0\) be a \(\pi\)-adically topological ring for an element \(\pi \in A_0\) and let \(A\) be a Tate ring \(A_0[1/\pi]\).
    Assume that \(A_0\) is integral over some Noetherian ring.
    Then \((A^{\widehat{u}})^\circ\) is the same as the \(\pi\)-adic completion \(\widehat{(A_0)_A^+}\) of the integral closure \((A_0)_A^+\) of \(A_0\) in \(A\).
\end{lemma}

\begin{proof}
    By \Cref{PowerBoundedUniformization}, we have \((A^{\widehat{u}})^\circ \cong \widehat{A^{u \circ}}\).
    Recall that \((A_0)_A^+\) becomes a ring of definition of \(A^u\) and thus \(A^{u \circ}\) is the complete integral closure \(((A_0)_A^+)_A^*\) of \((A_0)_A^+\) in \(A^u\) by \cite[Lemma 2.13 (1)]{nakazato2022Finite}.
    By assumption, \((A_0)_A^+\) is integral over some Noetherian ring and then \(A^{u \circ} = (A_0)_A^+\) by \cite[Proposition 7.1]{nakazato2023Variant}.
    This completes the proof.
\end{proof}

The main result of this paper is the following.

\begin{theorem} \label{Maintheorem2}
    Let \(R\) be a ring which contains \(\setZ\) as a subring and satisfies the following conditions:
    \begin{enumerate}
        \item The \(p\)-adic completion \(A_0 \defeq \widehat{R}\) of \(R\) has a map from some integral perfectoid ring.
        \item The perfectoidization \((\widehat{R})_{\perfd}\) of \(\widehat{R}\) is an (honest) integral perfectoid ring.
        \item The \(p\)-adic completion \(\widehat{C(R^{\ptf})}\) of the \(p\)-root closure \(C(R^{\ptf})\) is an integral perfectoid ring.
    \end{enumerate}

    Set a Tate ring \(A \defeq A_0[1/p]\).
    Then there exists an isomorphism \(\varphi \colon \widehat{C(R^{\ptf})} \xrightarrow{\cong} (\widehat{R})_{\perfd}^{\ptf}\) which is a restriction of the unique map \((A^{\widehat{u}})^\circ \to (A_0)_{\perfd, *}\) taken in \Cref{UniversalityUniform}.
    In particular, \((\widehat{R})_{\perfd}\) is \((p)^{1/p^\infty}\)-almost isomorphic to \(\widehat{C(R^{\ptf})}\).
    If \(R\) is \(p\)-torsion-free, we have an honest isomorphism \(\widehat{C(R)} \cong (\widehat{R})_{\perfd}\).



\end{theorem}

\begin{proof}
    Since \((A_0)_{\perfd}\) is an integral perfectoid ring whose perfectoid element is \(\varpi\),
    the ring \((A_0)_{\perfd}^{\ptf}\) is also an integral perfectoid ring by \cite[Proposition 2.19]{dospinescu2023Conjecture}.
    In particular, \((A_0)_{\perfd}^{\ptf}\) is \(p\)-root closed in \((A_0)_{\perfd}[1/p]\) because of the injectivity of the \(p\)-th power map \((A_0)_{\perfd}^{\ptf}/\varpi (A_0)_{\perfd}^{\ptf} \xrightarrow{a \mapsto a^p} (A_0)_{\perfd}^{\ptf}/p (A_0)_{\perfd}^{\ptf}\) (see \cite[(2.1.7.1)]{cesnavicius2023Purity}).
    The map of Tate rings \(R[1/p] \to (A_0)_{\perfd}[1/p]\) induced by \(R \to A_0 \to (A_0)_{\perfd}\) gives a unique map \(\map{\varphi}{\widehat{C(R^{\ptf})}}{(A_0)_{\perfd}^{\ptf}}\) such that the following diagram commutes
    \begin{center}
        \begin{equation} \label{RootClosurePerfectoidization}
            \begin{tikzcd}
                {R[1/p]} \arrow[r] & {A_0[1/p]} \arrow[r]   & {(A_0)_{\perfd}[1/p]}    \\
                & \widehat{C(R^{\ptf})} \arrow[rd, "\exists ! \varphi"] & \\
                C(R^{\ptf}) \arrow[uu, hook] \arrow[ru] \arrow[r] & C(A_0^{\ptf}) \arrow[r] \arrow[uu, hook, bend right]  & (A_0)_{\perfd}^{\ptf} \arrow[uu, hook] \\
                R \arrow[r] \arrow[u]     & A_0 \arrow[r] \arrow[u]    & (A_0)_{\perfd}. \arrow[u]
            \end{tikzcd}
        \end{equation}
    \end{center}

    Taking the \(p\)-adic completion of \(R \to C(R^{\ptf})\), we have a map \(A_0 \to \widehat{C(R^{\ptf})}\).
    By assumption, \(\widehat{C(R^{\ptf})}\) is an integral perfectoid ring and then, there exists a unique map \((A_0)_{\perfd} \to \widehat{C(R^{\ptf})}\) which extends \(A_0 \to \widehat{C(R^{\ptf})}\) as follows.
    \begin{center}
        \begin{equation} \label{PerfectoidizationRootClosure}
            \begin{tikzcd}
                (A_0)_{\perfd} \arrow[r, "\exists ! \psi"] & \widehat{C(R^{\ptf})}  \\
                A_0 \arrow[ru] \arrow[u]   &     \\
                R \arrow[r] \arrow[u]    & C(R^{\ptf}) \arrow[uu]
            \end{tikzcd}
        \end{equation}
    \end{center}
    Since \(\widehat{C(R^{\ptf})}\) is \(p\)-torsion-free, there exists a unique map \(\map{\psi}{(A_0)_{\perfd}^{\ptf})}{\widehat{C(R^{\ptf})}}\) which extends \((A_0)_{\perfd} \to \widehat{C(R^{\ptf})}\).
    
    Combining (\ref{RootClosurePerfectoidization}) and (\ref{PerfectoidizationRootClosure}), the composite map \((A_0)_{\perfd}^{\ptf} \xrightarrow{\psi} \widehat{C(R^{\ptf})} \xrightarrow{\varphi} (A_0)_{\perfd}^{\ptf}\) is an \(R\)-algebra map and thus an \(A_0 = \widehat{R}\)-algebra map.
    Furthermore, this extends to an \(A_0\)-algebra map of perfectoid Tate rings \((A_0)_{\perfd}[1/p] \to (A_0)_{\perfd}[1/p]\) and thus extends to an \(A_0\)-algebra map \((A_0)_{\perfd,*} \to (A_0)_{\perfd,*}\) by \Cref{UniformalityAlmostElement}.
    This must be the identity map because of \Cref{IdentityMap} above, and then \(\psi \circ \varphi\) is also the identity map. 
    
    On the other hand, consider the \(R\)-algebra map \(\widehat{C(R^{\ptf})} \xrightarrow{\varphi} (A_0)_{\perfd}^{\ptf} \xrightarrow{\psi} \widehat{C(R^{\ptf})}\).
    Inverting \(p\), we have a map of Tate rings \(\widehat{C(R^{\ptf})}[1/p] \to \widehat{C(R^{\ptf})}[1/p]\) over \(R[1/p]\) via \(R[1/p] \hookrightarrow C(R^{\ptf})[1/p] \to \widehat{C(R^{\ptf})}[1/p]\).
    Since \(C(R^{\ptf}) \subseteq R[1/p]\) is contained in the integral closure \(\widetilde{R}\) of \(R\) in \(R[1/p]\), the uniform completion \(C(R^{\ptf})[1/p]^{\widehat{u}}\) of \(C(R^{\ptf})[1/p]\) can be written as a Tate ring \(\widehat{\widetilde{R}}[1/p]\) by \Cref{UniformizationUniformCompletion}.
    Similarly, the uniform completion \(R[1/p]^{\widehat{u}}\) of \(R[1/p]\) also coincides with the Tate ring \(\widehat{\widetilde{R}}[1/p]\).
    So these two uniform completions \(C(R^{\ptf})[1/p]^{\widehat{u}}\) and \(R[1/p]^{\widehat{u}}\) are isomorphic each other via the canonical map of Tate rings \(R[1/p] \to C(R^{\ptf})[1/p]\).\footnote{Note that the identity map \(R[1/p] = C(R^{\ptf})[1/p]\) is only an isomorphism of rings which is not necessarily an isomorphism of topological rings. Therefore, we go back to the construction of uniform completion and prove the isomorphism \(R[1/p]^{\widehat{u}} \cong C(R^{\ptf})[1/p]^{\widehat{u}}\) in this way.}
    Furthermore, since \(\widehat{C(R^{\ptf})}\) is a \(p\)-torsion-free integral perfectoid ring, the uniform completion of \(\widehat{C(R^{\ptf})}[1/p]\) is isomorphic to itself by \Cref{PerfectoidRingTateRing}.
    These arguments show that
    \begin{equation}
        R[1/p]^{\widehat{u}} \xrightarrow{\cong} C(R^{\ptf})[1/p]^{\widehat{u}} \xrightarrow{\cong} \widehat{C(R^{\ptf})}[1/p]^{\widehat{u}} \cong \widehat{C(R^{\ptf})}[1/p]. \label{UniformCompRootClosure}
    \end{equation}

    This shows that the above \(R[1/p]\)-algebra map of Tate rings \(\widehat{C(R^{\ptf})}[1/p] \to \widehat{C(R^{\ptf})}[1/p]\) coincides with an \(R[1/p]\)-algebra map of Tate rings \(R[1/p]^{\widehat{u}} \to R[1/p]^{\widehat{u}}\) and this is the identity map by the universality of uniform completion.

    We next check that \(\map{\varphi}{\widehat{C(R^{\ptf})}}{(A_0)_{\perfd}^{\ptf}}\) is a restriction of the unique map \((A^{\widehat{u}})^\circ \to (A_0)_{\perfd,*}\) taken in \Cref{UniversalityUniform}.
    Similarly as above, the uniform completion of \(A_0[1/p]\) is
    \begin{equation}
        A^{\widehat{u}} = A_0[1/p]^{\widehat{u}} = \widehat{R}[1/p]^{\widehat{u}} \cong R[1/p]^{\widehat{u}} \cong \widehat{C(R^{\ptf})}[1/p]^{\widehat{u}}.
    \end{equation}
    The uniform completion \(\widehat{C(R^{\ptf})}[1/p]^{\widehat{u}}\) of \(\widehat{C(R^{\ptf})}[1/p]\) has a universality which gives a unique extension \(A^{\widehat{u}} \to (A_0)_{\perfd}[1/p]\) of \(\varphi\) as follows (see \Cref{UniversalityUniformCompletion}):
    \begin{center}
        \begin{tikzcd}
            R \arrow[r] \arrow[d]      & A_0 \arrow[r]   & (A_0)_{\perfd}^{\ptf} \arrow[d, hook] \\
            \widehat{C(R^{\ptf})} \arrow[rru, "\varphi"] \arrow[d, hook]         &   & {(A_0)_{\perfd}[1/p]}   \\
            {\widehat{C(R^{\ptf})}[1/p]} \arrow[d] \arrow[rru, "{\varphi[1/p]}"] & &       \\
            {\widehat{C(R^{\ptf})}[1/p]^{\widehat{u}}} \arrow[r, "\cong"]      & A^{\widehat{u}}. \arrow[ruu, "\exists ! \eta", dashed] &                         
        \end{tikzcd}
    \end{center}
    Since \(A_0\) is the \(p\)-adic completion of \(R\), \(\map{\eta}{A^{\widehat{u}}}{(A_0)_{\perfd}[1/p]}\) is a unique extension of \(A_0 \to (A_0)_{\perfd}^{\ptf}\).
    By \Cref{UniformalityAlmostElement} and the proof of \Cref{UniversalityUniform}, \(\eta\) induces the unique map \((A^{\widehat{u}})^\circ \to (A_0)_{\perfd,*}\) taken in \Cref{UniversalityUniform} and this extends \(\varphi\).

    If \(R\) is \(p\)-torsion-free, \(A_0 = \widehat{R}\) is also \(p\)-torsion-free and so is \((A_0)_{\perfd}\) by \cite[Lemma A.2]{ma2022Analogue}.
    Then \(\widehat{C(R^{\ptf})} = \widehat{C(R)}\) is isomorphic to \((\widehat{R})_{\perfd}\).
\end{proof}




\begin{remark} \label{pRootClosurePerfectoid}
    Let \(R\) be a (not necessarily \(p\)-adically complete) ring containing a compatible sequence of \(p\)-power roots \(\{\varpi^{1/p^j}\}_{j \geq 0}\) of \(\varpi \in R\) such that \(\varpi^p\) is some unit multiple of \(p\) in \(R\).
    If the Frobenius map \(R/pR \xrightarrow{F} R/pR\) is surjective, the \(p\)-adic completion \(\widehat{C(R)}\) is an integral perfectoid ring by \cite[Proposition 2.1.8]{cesnavicius2023Purity}.
    Compare \Cref{UniformCompletionPerfectoid} above.
    
    In particular, any semiperfectoid ring satisfies the assumptions of \Cref{Maintheorem2} and so we have the following corollary.
\end{remark}

\begin{corollary} \label{MainTheorem3}
    Let \(R\) be a ring such that the \(p\)-adic completion \(\widehat{R}\) of \(R\) becomes a semiperfectoid ring which contains \(\setZ\) as a subring.
    Then \((\widehat{R})_{\perfd}^{\ptf}\) is isomorphic to the \(p\)-adic completion \(\widehat{C(R^{\ptf})}\) of \(C(R^{\ptf})\).
    If \(R\) is \(p\)-torsion-free, we have \((\widehat{R})_{\perfd} \cong \widehat{C(R)}\).

\end{corollary}

\begin{proof}
    First, the semiperfectoid ring \(\widehat{R}\) has a surjective map from some integral perfectoid ring by the definition of semiperfectoid rings (\Cref{Semiperfectoid}).
    Second, The perfectoidization \((\widehat{R})_{\perfd}\) of \(\widehat{R}\) is an (honest) integral perfectoid ring by \cite[Corollary 7.3 and Proposition 8.5]{bhatt2022Prismsa}.
    Finally, the \(p\)-adic completion \(\widehat{C(R^{\ptf})}\) of the \(p\)-root closure \(C(R^{\ptf})\) is an integral perfectoid ring by \Cref{pRootClosurePerfectoid}.
    So \(R\) satisfies all assumptions of \Cref{Maintheorem2}, and this completes the proof.


\end{proof}

\section{Connections between \texorpdfstring{\(p\)}{p}-root Closure and Perfectoidization}

We recall a mixed characteristic analog of the perfection of rings, which was introduced in \cite{ishizuka2023Mixed}.
This construction includes an example from \cite[Section 4]{roberts2008Root} which demonstrates a good behavior of \(p\)-root closure from the perspective of Fontaine rings.

\begin{construction}
    \label{ExplicitConst}
    Let \((R_0,\mfrakm,k)\) be a complete Noetherian local domain of mixed characteristic \((0, p)\) with perfect residue field \(k\) and let \(p, x_2, \dots, x_n\) be any system of generators of the maximal ideal \(\mfrakm\) such that \(p, x_2, \dots , x_d\) forms a system of parameters of \(R_0\).
    Choose compatible sequences of \(p\)-power roots 
    \begin{equation}
        \{p^{1/p^j}\}_{j \ge 0},\{x_2^{1/p^j}\}_{j \ge 0},\ldots,\{x_n^{1/p^j}\}_{j \ge 0}
    \end{equation}
    inside the absolute integral closure \(R_0^+\).
    Set
    \begin{equation}
        R_{\infty} := \bigcup_{j \ge 0} R_0[p^{1/p^j},x_2^{1/p^j},\ldots,x_n^{1/p^j}] \subseteq R_0^+.
    \end{equation}
    Let \(\widetilde{R}_{\infty}\) be the integral closure of \(R_{\infty}\) in \(R_{\infty}[1/p]\) and let \(\widehat{\widetilde{R}}_{\infty}\) (resp, \(\widehat{R}_{\infty}\)) be the \(p\)-adic completion of \(\widetilde{R}_{\infty}\) (resp, \(R_{\infty}\)).

    By Cohen's structure theorem, there exists a surjective map \(S_0 \twoheadrightarrow R_0\) such that \(S_0\) is a complete unramified regular local ring \(W(k)[|t_2, \dots, t_n|]\) and \(t_i\) maps to \(x_i\) respectively.
    In particular, \(R_0\) is a finite extension of a subring \(T_0 \defeq W(k)[|t_2, \dots, t_d|]\) of \(S_0\).
    Then, proceeding as in \cite{ishizuka2023Mixed}, we have a surjective map \(S_{\infty} \twoheadrightarrow R_{\infty}\) and its \(p\)-adic completion \(\widehat{S}_{\infty} \to \widehat{R}_{\infty}\).
    Remark that \(\widehat{S}_{\infty}\) is an integral perfectoid ring which has a compatible sequence of \(p\)-power roots \(\{p^{1/p^j}\}_{j \geq 0}\) of \(p\).
    Therefore, \(\widehat{R}_{\infty}\) is a semiperfectoid ring.

    In \cite{ishizuka2023Mixed}, we show that \(\widehat{C(R_{\infty})}\) and \(\widehat{\widetilde{R}}_{\infty}\) are integral perfectoid rings and are \((pg)^{1/p^\infty}\)-almost flat and \((pg)^{1/p^\infty}\)-almost faithful \(T_0\)-algebra where \(g\) is a non-zero element of \(\widehat{R}_{\infty}\) that becomes a non-zero-divisor in \(\widehat{\widetilde{R}}_{\infty}\).
\end{construction}

\begin{corollary} \label{ExplicitPerfectoidization}
    Keep the notation of \Cref{ExplicitConst}.
    Then, the integral perfectoid ring \(\widehat{\widetilde{R}}_{\infty}\) is \((p)^{1/p^\infty}\)-almost isomorphic to the perfectoidization \((\widehat{R}_{\infty})_{\perfd}\) and its restriction to \(\widehat{C(R_{\infty})}\) induces an isomorphism of \(\widehat{C(R_{\infty})} \xrightarrow{\cong} (\widehat{R}_{\infty})_{\perfd}\).
    
\end{corollary}

\begin{proof}
    As above, \(\widehat{R}_{\infty}\) is a semiperfectoid ring and has a compatible sequence of \(p\)-power roots \(\{p^{1/p^j}\}_{j \geq 0}\) of \(p\).
    So the second statement is already proved in \Cref{MainTheorem3}.
    
    In particular, the perfectoidization \((\widehat{R}_{\infty})_{\perfd}\) is an honest integral perfectoid ring by \cite[Proposition 8.5]{bhatt2022Prismsa}.
    Also \(\widehat{\widetilde{R}}_{\infty}\) is an integral perfectoid ring by \cite[Lemma 4.2]{ishizuka2023Mixed}.
    Note that the proof of \cite[Proposition 5.9]{ishizuka2023Mixed} shows that the uniform completion of \(\widehat{R}_{\infty}[1/p]\) is isomorphic to the uniform completion \(R_{\infty}[1/p]^{\widehat{u}}\) of \(R_{\infty}[1/p]\).
    Since \(\widehat{R}_{\infty}\) is a semiperfectoid ring, the uniform completion \(\widehat{R}_{\infty}[1/p]^{\widehat{u}} \cong R_{\infty}[1/p]^{\widehat{u}}\) is a perfectoid Tate ring by \Cref{UniformCompletionPerfectoid}.
    In particular, we have \((R_{\infty}[1/p]^{\widehat{u}})^\circ = \widehat{\widetilde{R}}_{\infty}\) by \Cref{UniformCompRealization} and deduce the first statement by \Cref{MainTheorem1}.
    %
\end{proof}



\end{document}